%% file: cvx_review.tex
\documentclass[sts, preprint]{imsart}
\setattribute{journal}{name}{~}

\usepackage{amsthm,amsmath,natbib,color,prettyref, bm, amsfonts, amssymb, cases}
\usepackage[breaklinks=true]{hyperref}
\usepackage{breakcites}
\RequirePackage[colorlinks,citecolor=blue,urlcolor=blue]{hyperref}

\usepackage{xspace}
\startlocaldefs
\numberwithin{equation}{section}
\theoremstyle{plain}
\newtheorem{theorem}{Theorem}[section]

\newtheorem{lemma}{Lemma}[section]

\theoremstyle{definition}

\newtheorem{remark}{Remark}[section]
\endlocaldefs

\input{macros}

\begin{document}

\begin{frontmatter}

\title{Convex Relaxation Methods for Community Detection}
\runtitle{Convex Community Detection}

\begin{aug}
	\author{\fnms{Xiaodong} \snm{Li}\ead[label=e1]{xdgli@ucdavis.edu}},
	\author{\fnms{Yudong} \snm{Chen}\ead[label=e2]{yudong.chen@cornell.edu}}
	\and
	\author{\fnms{Jiaming} \snm{Xu}\ead[label=e3]{jx77@duke.edu}}
	
	\runauthor{X. Li et al.}
	
	\affiliation{University of California, Davis,~Cornell University~and~Duke University}
	
	\address{Xiaodong Li is an assistant professor, Department of Statistics, UC Davis, 399 Crocker Lane, Math Sciences Bldg, Davis, CA, 95616  \printead{e1}.}
	\address{Yudong Chen is an assistant professor, School of Operations Research and Information Engineering, Cornell University, 223 Frank H.T. Rhodes Hall, Ithaca, NY 14853  \printead{e2}.}	
	\address{Jiaming Xu is an assistant professor, Fuqua School of Business, Duke University, 100 Fuqua Drive, Durham, NC 27708, \printead{e3}.}
	
\end{aug}



\runauthor{X. Li et al.}

\begin{abstract}
	This paper surveys recent theoretical advances in convex optimization approaches for community detection. We introduce some important theoretical techniques and results for establishing the consistency of convex community detection under various statistical models. In particular, we discuss the basic techniques based on the primal and dual analysis.
    We also present results that demonstrate several distinctive advantages of convex community detection, including robustness against outlier nodes, consistency under weak assortativity, and adaptivity to heterogeneous degrees. 
    
    This survey is not intended to be a complete overview of the vast literature on this fast-growing topic. Instead, we aim to  provide a big picture of the remarkable recent development in this area and to make the survey accessible to a broad audience. We hope that this expository article can serve as an introductory guide for readers who are interested in using, designing, and analyzing convex relaxation methods in network analysis. 
\end{abstract}



\end{frontmatter}

\input{intro}

\input{dual2}

\input{sharp}

\input{weak_assortativity}

\input{primal}

\input{cmm}

\input{robustness}

\input{conclusion}

%
%
%
%

\bibliographystyle{imsart-nameyear}
\bibliography{refs}

\end{document}

%% file: macros.tex
\newcommand{\Perp}{\indep}

\newcommand{\reals}{{\mathbb{R}}}

\newcommand{\expect}[1]{\mathbb{E}\left[ #1 \right]}

\newcommand{\prob}[1]{{ \mathbb{P}\left\{ #1 \right\} }}

\newcommand\indep{\protect\mathpalette{\protect\independenT}{\perp}}
\def\independenT#1#2{\mathrel{\rlap{$#1#2$}\mkern2mu{#1#2}}}

\newcommand{\Binom}{{\rm Binom}}

\newrefformat{eq}{(\ref{#1})}
\newrefformat{chap}{Chapter~\ref{#1}}
\newrefformat{sec}{Section~\ref{#1}}
\newrefformat{algo}{Algorithm~\ref{#1}}
\newrefformat{fig}{Fig.~\ref{#1}}
\newrefformat{tab}{Table~\ref{#1}}
\newrefformat{rmk}{Remark~\ref{#1}}
\newrefformat{clm}{Claim~\ref{#1}}
\newrefformat{def}{Definition~\ref{#1}}
\newrefformat{cor}{Corollary~\ref{#1}}
\newrefformat{lmm}{Lemma~\ref{#1}}
\newrefformat{prop}{Proposition~\ref{#1}}
\newrefformat{app}{Appendix~\ref{#1}}
\newrefformat{hyp}{Hypothesis~\ref{#1}}
\newrefformat{thm}{Theorem~\ref{#1}}

\newcommand{\sth}[1]{\left\{ #1 \right\}}

\newcommand{\norm}[1]{\left\|{#1} \right\|}

\newcommand{\iprod}[2]{\left \langle #1, #2 \right\rangle}
\newcommand{\Iprod}[2]{\langle #1, #2 \rangle}

\newcommand{\calZ}{{\mathcal{Z}}}

\newcommand{\diag}{{\rm diag}}

\newcommand{\allones}{\mtx{J}}
\newcommand{\identity}{\mtx{I}}

\definecolor{darkred}{RGB}{100,0,0}
\definecolor{darkgreen}{RGB}{0,100,0}
\definecolor{darkblue}{RGB}{0,0,150}
\definecolor{red}{RGB}{255,0,0}

\newcommand{\E}{\operatorname{\mathbb{E}}}

\renewcommand{\vec}[1]{{\boldsymbol{#1}}}

\newcommand{\bM}{\bm{M}}

\newcommand{\vct}[1]{\bm{#1}}
\newcommand{\mtx}[1]{\bm{#1}}

\newcommand{\bee}{\begin{equation*}}
\newcommand{\eee}{\end{equation*}}
\newcommand{\bea}{\begin{eqnarray}}
\newcommand{\eea}{\end{eqnarray}}
\newcommand{\beaa}{\begin{eqnarray*}}
\newcommand{\eeaa}{\end{eqnarray*}}

\renewcommand{\hat}{\widehat}

\newcommand{\SBM}{{\sf SBM}\xspace}
\newcommand{\DCSBM}{{\sf DCSBM}\xspace}

\newcommand{\silent}[1]{}

\definecolor{xl}{RGB}{200,50,120}
\definecolor{yc}{RGB}{50,150,50}
\definecolor{jx}{RGB}{50,200,200}

%% file: intro.tex

\section{Introduction}
\label{sec:intro}
Convex relaxation has arisen as a powerful framework for developing computationally tractable and statistically efficient solutions to the community detection problems. Particularly in the last few years, this area has enjoyed remarkable progress: a variety of new methods have been proposed with established strong performance guarantees demonstrating their power and statistical advantages. In this expository article, we give a survey of these convex optimization approaches from the perspective of convexified maximum likelihood and discuss their major theoretical properties as well as relevant analytical tools. We hope this exposition is helpful to the readers who are interested in proposing and analyzing their own convex optimization methods for network analysis.


Convex optimization approaches for community detection can be traced back to the computer science and mathematical programming literature in the study of the planted partition problem; see, e.g., \citet{MS2010, OH2011,Jalali2011clustering, AV2011, Ames2013}. For community detection under statistical models, various theoretical properties of convex optimization methods have been studied in depth recently; a partial list includes strong consistency\footnote{An estimator of the community structure is said to achieve \emph{strong consistency}, if it equals the underlying
true community partition with probability tending to $1$ as network size grows. It is said to achieve \emph{weak consistency}, if the fraction of misclassified nodes (up to a permutation of community labels) goes to $0$ with probability tending to $1$. It is said to achieve \emph{non-trivial recovery}, if 
the fraction of misclassified nodes (up to a permutation of community labels) is strictly better than 
the random guess.} with a growing number of communities \citep{CSX2014,ChenXu14,Cai2014robust}, sharp threshold under sparse networks for strong consistency \citep{Abbe14, HajekWuXuSDP14, Bandeira15, HajekWuXuSDP15, ABBK, perry2015semidefinite}, weak consistency \citep{Vershynin14, 
FeiChen2017_exponential}, non-trivial recovery \citep{MontanariSen15,FanMontanari16,Mei17}, 
robustness against outlier nodes \citep{Cai2014robust,moitra2016robust,makarychev2016learning}, consistency under degree-corrected models \citep{CLX2015}, and consistency under weak assortativity \citep{AL2014, YSC2018}.


\subsection{Stochastic Block Models}

In the above literature, a standard setting for deriving and analyzing the convex relaxation formulations is the so-called \emph{stochastic block model} (\SBM) \citep{HLL1983}, which we shall focus on. The \SBM~is a generative model for a random graph $G=(V, E)$, where $V=[n]$ is a set of $n$ nodes and and $E$ is the set of edges. Under \SBM, the nodes are partitioned into $r$ clusters according to the mapping $\phi: [n] \rightarrow [r]$, where the node $i$ belongs to the cluster $\phi(i)$. Each pair of nodes $i$ and $j$ are connected independently with probability $B_{\phi(i) \phi(j)}$. Note that the connectivity probability only depends on the cluster memberships of the nodes. We refer to the symmetric matrix $\mtx{B}:=[B_{ab}]_{1 \leq a, b \leq r}$ as the \emph{connectivity probability matrix}. We also denote the size of the $a$th cluster by $n_a:= |\phi^{-1}(a)|$. The general community detection problem under \SBM~is to estimate the unknown cluster structure $ \phi $ given the observed graph $ G $.

\SBM~is a powerful and versatile modeling tool that captures several key aspects of the community detection problem in real world networks:
\begin{itemize}
\item \emph{Sparsity}: Most large-scale networks are sparse. In \SBM, sparsity can be modeled by assuming that the  matrix $\mtx{B}$ has small entries.
\item \emph{Connectivity probabilities}: Different connectivity patterns can be captured by the matrix $\mtx{B}$. Examples include \emph{strong assortativity}, where $\min_{1 \leq a \leq r} B_{aa} > \max_{1 \leq  a < b \leq r} B_{ab}$, and \emph{weak assortativity} $B_{aa} - \max_{b \neq a} B_{ab} >0, \forall a$ \citep{AL2014}.
\item \emph{Number of communities}: The number of clusters $r$ may be large and in particular is allowed to grow with the number of nodes $ n $.
\item \emph{Unbalancedness}: \SBM~allows for significantly unbalanced clusters by imposing high variation in the cluster sizes $n_1, \ldots, n_r$.
\end{itemize}

For expository purposes, we mostly focus on the simple $(p, q)$-\SBM~setting, where the diagonal entries of $\mtx{B}$ are all $p$ and the off-diagonal entries are all $q$. In other words, nodes in the same clusters are connected with probability $p$, whereas nodes in different clusters are connected with probability $q$. The results and techniques discussed in this paper can often be extended to the general \SBM~(or at least to the strong assortativity setting), though in some cases the extension is more difficult.

To prepare for discussion of the convex relaxation approaches, we introduce an alternative way to parameterize the cluster structure through an $ n\times n $ \emph{partition matrix} $\mtx{X}$. 
Here the binary variable $X_{ij}$ indicates whether or not the nodes $i$ and $j$ are assigned to the same cluster, i.e. $X_{ij}=1$ if $\phi(i) = \phi(j) $ and $X_{ij}=0$ if $\phi(i) \neq \phi(j)$. Correspondingly, the observed network can also be represented by its $ n\times n $ adjacent matrix $ \mtx{A} $, where $A_{ij}=1$ if the nodes $i$ and $j$ are connected and $ A_{ij} = 0 $ otherwise.

\subsection{Convex Relaxation Approaches}

Convex relaxations for community detection can be derived in various ways, for instance through the relaxation of 
modularity maximization, $k$-means, or min $k$-cut. Such derivations often give rise to convex programs of similar forms. Here we focus on the perspective of convexifying maximum likelihood estimators. To the best of our knowledge, this perspective was first considered in the context of \SBM~by \citet{CSX2014}.

Let us start with the $(p,q)$-\SBM, and derive the likelihood function $ \ell(\mtx{A}|\mtx{X}, p, q)$ of the community structure given the observed network data. If we denote by $p(A_{ij} |X_{ij}, p, q)$ the probability mass function for the Bernoulli random variable $A_{ij}$, then by definition of the model we have: 
\begin{align*}
\log p(A_{ij}|X_{ij}, p, q) &= \left\{\begin{aligned}\log p&,\quad A_{ij}=X_{ij}=1 \\ \log (1-p)&, \quad A_{ij}=0, X_{ij}=1 \\ \log q&,\quad A_{ij}=1, X_{ij}=0 \\ \log (1-q)&, \quad A_{ij}=X_{ij}=0. \end{aligned} \right.
\end{align*}
This formula can be rewritten as a linear function of $ X_{ij} $:
\begin{align*}
\log p(A_{ij}|X_{ij}, p, q) 
&= \left[\left(\log \frac{p}{1-p} - \log \frac{q}{1-q}\right)A_{ij} - \log \frac{1-q}{1-p}\right]X_{ij} 
\\
&~~+ \left[A_{ij} \log q + (1-A_{ij})\log(1-q)\right].
\end{align*}
Since all pairs of nodes are connected independently, that is, the $A_{ij}$'s are independent across $1 \leq i < j \leq n$, the log-likelihood function given the observed adjacency matrix $ \mtx{A} $ has the summation form
\[
\ell(\mtx{A}|\mtx{X}, p, q) = \sum_{1 \leq i < j \leq n} \log p(A_{ij}|X_{ij}, p, q).
\]

Our goal is to estimate the true matrix of membership, denoted by $\mtx{X}^\star$, using the maximum likelihood approach.
If $p$ and $q$ are fixed and given, then maximizing the log likelihood $\ell(\mtx{A}|\mtx{X}, p, q)$ is equivalent to the maximization of 
\[
\sum_{1 \leq i < j \leq n} \left[\left(\log \frac{p}{1-p} - \log \frac{q}{1-q}\right)A_{ij} - \log \frac{1-q}{1-p}\right]X_{ij} 
\]
over all possible cluster matrices $\mtx{X}$ with any number of clusters $r$.
In real world networks, individuals within the same communities are often more likely to be connected than those across different communities, which corresponds to the setting with $p>q$ and hence $\log \frac{p}{1-p} - \log \frac{q}{1-q}>0$. In this case, the above optimization problem is equivalent to 
\begin{equation}
\label{eq:likelihood}
\max_{\mtx{X}} \sum_{1 \leq i < j \leq n}  (A_{ij} - \lambda) X_{ij} \; ,
\end{equation}
where
\begin{equation}
\label{eq:lambda_p_q}
\lambda \equiv \lambda(p,q) =\frac{\log(1-q)-\log (1-p)}{\log p-\log q+\log (1-q)-\log (1-p)}.
\end{equation}
If $p$ and $q$ are unknown, we simply treat $\lambda$ as a tuning parameter.\footnote{An alternative approach is profile likelihood maximization, i.e., we maximize $\ell(\mtx{A}|\mtx{X}, \hat{p}, \hat{q})$ with $(\hat{p}(\mtx{X}), \hat{q}(\mtx{X})):= \arg\max_{p, q} \ell(\mtx{A}|\mtx{X}, p, q)$. However, the framework of profile likelihood will result in a highly nonlinear function of $\mtx{X}$ to maximize. }

To derive a convex relaxation, it is more convenient to write the problem~\prettyref{eq:likelihood} in a matrix form. Since both $\mtx{A}$ and $\mtx{X}$ are symmetric and the diagonal entries of $\mtx{A}$ are all $0$, the problem \prettyref{eq:likelihood} can be rewritten as
\begin{equation}
\label{eq:likelihood_matrix}
\begin{aligned}
&\max_{\mtx{X}} &&\left\langle \mtx{X}, \mtx{A}-\lambda \mtx{J}_n \right\rangle
\\
&\text{subject to} && \mtx{X} \text{~is a partition matrix},
\end{aligned}
\end{equation}
where $\mtx{J}_n$ is the $n \times n$ matrix with all entries equal to $1$, and $ \langle \mtx{A}, \mtx{B} \rangle := \text{trace}(\mtx{A}^\top \mtx{B}) $ denotes the trace inner product.
The problem \prettyref{eq:likelihood_matrix} is non-convex due to its constraint. To convexify this program, let us investigate the properties of a community matrix $\mtx{X}$. By definition we see that $ \mtx{X} $ must have form
\begin{equation}
\label{eq:membership}
\mtx{X}=\mtx{P}\begin{bmatrix} \mtx{J}_{n_1} &~ & ~\\ ~ &\ddots & ~\\ ~ & ~ & \mtx{J}_{n_r}\end{bmatrix}\mtx{P}^{\intercal},
\end{equation}
where $\mtx{P}$ is some permutation matrix and the number of communities $r$ is unknown. The set of all matrices $ \mtx{X} $ of this form is of course non-convex. The key observation is that any such $\mtx{X}$ satisfies several convex constraints such as
(i) all entries of $\mtx{X}$ are nonnegative,
(ii) all diagonal entries of $\mtx{X}$ are $1$, and
(iii) $\mtx{X}$ is positive semi-definite.
Various convex relaxations can hence be derived by replacing the constraint in~\prettyref{eq:likelihood_matrix} with these convex constraints (or variants thereof). 

One typical relaxation obtained in this way is:
\begin{equation}
\label{eq:SDP1}
\begin{aligned}
&\text{max} && \left\langle \mtx{X}, \mtx{A}-\lambda \mtx{J}_n \right\rangle
\\
&\text{subject to} && \mtx{X} \succeq \mtx{0}, ~\mtx{X} \geq \mtx{0}, ~X_{ii} =1 \text{~for~} 1\leq i \leq n,
\end{aligned}
\end{equation}
where $\mtx{X} \geq \mtx{0}$ means $X_{ij} \geq 0$ for $1\leq i, j\leq n$. 
In \prettyref{sec:dual} we present results on when \prettyref{eq:SDP1} recovers the ground truth clusters. These results highlight the following attractive properties of the formulation \prettyref{eq:SDP1}:
\begin{itemize}
\item The communities are allowed to be significantly unbalanced;
\item The number of communities $ r $ may grow as $n$ increases;
\item Although \prettyref{eq:SDP1} is derived from the $(p, q)$-model, it is applicable to a more general connectivity probability matrix $\mtx{B}$ with strong assortativity;
\item The knowledge of $r$ is not required in \prettyref{eq:SDP1};
\item There is only one tuning parameter, $\lambda$.
\end{itemize}

Convex relaxation methods of similar forms to \prettyref{eq:SDP1} have been extensively studied in the literature; see, e.g., \citet{Ames2013, ChenXu14, Cai2014robust, AL2014, Vershynin14, HajekWuXuSDP15,  ABBK, perry2015semidefinite, CLX2015, FeiChen2017_exponential}. Moreover, other convex relaxation formulations, such as those for the simple two-community setting considered by \citet{Abbe14, HajekWuXuSDP14, Bandeira15, MontanariSen15, JavamardMotanariRicci15}, are also intrinsically related to \prettyref{eq:SDP1}. 

The subsequent sections are devoted to understanding the statistical performance of  the convex relaxations~\prettyref{eq:SDP1} and some of its variants. In particular, we present theoretical results on when these relaxations yield, with high probability, a good estimate of the ground truth communities of \SBM~in terms of strong and weak consistency. Particular focus is given to elucidating several important techniques, including dual and primal analysis, for establishing such results.

In the sequel, we use $ C_0 $ etc.\ to denote a positive and sufficiently large absolute constant. By ``with high probability'' we mean with probability at least $ 1-10n^{-10} $.


%% file: dual2.tex

\section{Strong Consistency via Dual Analysis}
\label{sec:dual}
In this section we will introduce a generic dual analysis technique, based on the KKT condition, that can be used to establish strong consistency of the convex relaxation approach in \prettyref{eq:SDP1}. This type of arguments are widely employed in the literature to study exact recovery under \SBM~using convex optimization; see, e.g., \citet{CSX2014,Cai2014robust,Ames2013}. 
%
%

For the ease of presentation, we introduce the dual analysis under the two-cluster $(p, q)$-model with $n_1 = O(n)$ and $n_2 = O(n)$. Without loss of generality we assume that the permutation matrix corresponding to the ground truth communities is the identity, i.e., $\mtx{P} = \mtx{I}$. Under this assumption, the adjacency matrix $\mtx{A}$ and the true partition matrix $\mtx{X}^\star$ as defined in \prettyref{eq:membership} can be written as 
\begin{equation}
\label{eq:primal_two}
\mtx{A}=
\begin{bmatrix} 
\mtx{A}_{11} & \mtx{A}_{12}
\\
\mtx{A}_{12}^\top & \mtx{A}_{22}
\end{bmatrix}
\text{~~and~~}
\mtx{X}^\star= \begin{bmatrix}
	\mtx{J}_{n_1} & \mtx{0} \\
	\mtx{0} & \mtx{J}_{n_2}
\end{bmatrix},
\end{equation}
where the entries of $\mtx{A}_{11}$ and $\mtx{A}_{22}$ are Bernoulli random variables with parameter $p$ and those $\mtx{A}_{12}$ are Bernoulli with parameter $q$.

Our goal is to prove that the convex relaxation~\prettyref{eq:SDP1} recovers the true $ \mtx{X}^\star $ as an optimal solution. We do so by showing that $\mtx{X}^\star $ satisfies the KKT condition (cf. \citealt[Section 5.5.3]{boyd2004}), a sufficient condition for optimality. The KKT condition stipulates the existence of some $\mtx{\Lambda},$ $\mtx{\Phi},$ and $\vct{\beta}$ satisfying 
\begin{itemize}
    \item Stationarity: $\lambda \mtx{J}_n - \mtx{A} = \mtx{\Lambda} + \mtx{\Phi} + \mbox{diag}(\vct{\beta})$,
	\item Complementary slackness: $\langle \mtx{X}^\star, \mtx{\Lambda}\rangle=\langle \mtx{X}^\star, \mtx{\Phi}\rangle=0$,
    \item Dual feasibility: $\mtx{\Phi} \ge \mtx{0}, \mtx{\Lambda} \succeq \mtx{0}$, 
\end{itemize}
where $\mtx{\Lambda}, \mtx{\Phi} \in \mathbb{R}^{n\times n}$ and $\vct{\beta} \in \mathbb{R}^n$ are the dual variables associated with
the constraints of \prettyref{eq:SDP1} in order. 
We write these variables in blocks as in \prettyref{eq:primal_two}:
\[
\mtx{\Lambda} =\begin{bmatrix}
	\mtx{\Lambda}_{11} & \mtx{\Lambda}_{12} \\
	\mtx{\Lambda}_{12}^\top & \mtx{\Lambda}_{22}
\end{bmatrix}, \quad \mtx{\Phi }=\begin{bmatrix}
	\mtx{\Phi}_{11} & \mtx{\Phi}_{12} \\
	\mtx{\Phi}_{12}^\top & \mtx{\Phi}_{22}
	\end{bmatrix}, \text{~~and~~}\vct{\beta} =\begin{bmatrix}
	\vct{\beta_{1}}  \\
	\vct{\beta_{2}}
\end{bmatrix}.
\]
We proceed by simplifying the KKT condition. The complementary slackness and dual constraints imply the following equalities: 
\[
\mtx{\Lambda}_{11}\vct{1}_{n_1}= \vct{0}, 
~~\mtx{\Lambda}_{22}\vct{1}_{n_2}= \vct{0},
~~\mtx{\Lambda}_{12}\vct{1}_{n_2}=\vct{0},
~~\mtx{\Lambda}_{12}^\top \vct{1}_{n_1}=\vct{0},
~~\mtx{\Phi}_{11}= \mtx{0},~\mtx{\Phi}_{22}=\mtx{0}.
\]
The stationary condition can be similarly rewritten in block structures
\[
	\left\{ \begin{aligned}
	&-\mtx{A}_{11} + \lambda\mtx{J}_{n_1} -\mtx{\Lambda}_{11}-\mbox{diag}(\vct{\beta}_1)=0, \\
	&-\mtx{A}_{22} + \lambda\mtx{J}_{n_2} -\mtx{\Lambda}_{22}-\mbox{diag}(\vct{\beta}_2)=0, \\
	&-\mtx{A}_{12} + \lambda\mtx{J}_{n_1, n_2} -\mtx{\Lambda}_{12}-\mtx{\Phi}_{12}=0.
	\end{aligned} \right.
\] 
Combining with $\mtx{\Lambda}_{11}\vct{1}_{n_1}= \vct{0}$ and $\mtx{\Lambda}_{22}\vct{1}_{n_2}= \vct{0}$, we can express $\mtx{\Lambda}$ by $\mtx{\Phi}_{12}$ and $\lambda$ as 
\[
\mtx{\Lambda} =\begin{bmatrix}
		\mbox{diag}(\mtx{A}_{11}\vct{1}_{n_1})-n_1\lambda \mtx{I}_{n_1}-\mtx{A}_{11}+\lambda \mtx{J}_{n_1}  & -\mtx{\Phi}_{12}-\mtx{A}_{12}+\lambda \mtx{J}_{n_1,n_2} 
        \\
		-\mtx{\Phi}_{12}^\top - \mtx{A}_{12}^\top + \lambda \mtx{J}_{n_2,n_1} & \mbox{diag}(\mtx{A}_{22}\vct{1}_{n_2})-n_2\lambda \mtx{I}_{n_2}-\mtx{A}_{22}+\lambda \mtx{J}_{n_2} 
	\end{bmatrix}.
\]	
By the conditions $\mtx{\Lambda}_{12}\vct{1}_{n_2}=\vct{0}$ and $\mtx{\Lambda}_{12}^\top \vct{1}_{n_1}=\vct{0}$, we have
\begin{equation}
\label{eq:kkt-cd}
\mtx{\Phi}_{12}\vct{1}_{n_2}=(\lambda\mtx{J}_{n_1, n_2}-\mtx{A}_{12})\vct{1}_{n_2}, \quad
\vct{1_}{n_1}^\top \mtx{\Phi}_{12}=\vct{1}_{n_1}^\top (\lambda\mtx{J}_{n_1, n_2}-\mtx{A}_{12}).
\end{equation}
Therefore, to establish the KKT condition, it suffices to construct a nonnegative matrix $\mtx{\Phi}_{12}\geq \mtx{0}$ such that $\mtx{\Lambda} \succeq \mtx{0}$ and \prettyref{eq:kkt-cd} holds. 

\paragraph{Construction of $\mtx{\Phi}_{12}$:} 
Such construction is often an art. Here we propose as a candidate using the following rank-$2$ matrix  
$\mtx{\Phi}_{12}=\vct{y} \vct{1}^\top_{n_2} + \vct{1}_{n_1} \vct{z}^\top$ for some vectors $\vct{y} \in \reals^{n_1}$
and $\vct{z}\in \reals^{n_2}$. We solve for $\vct{y}$ and $\vct{z}$ by plugging this candidate into \prettyref{eq:kkt-cd}, which yields that
\begin{equation*}
\begin{aligned}
\mtx{\Phi}_{12}&=-\frac{1}{n_2}\mtx{A}_{12}\mtx{J}_{n_2}-\frac{1}{n_1}\mtx{J}_{n_1} \mtx{A}_{12}+\frac{1}{n_1n_2}\mtx{J}_{n_1} \mtx{A}_{12}\mtx{J}_{n_2}+\lambda \mtx{J}_{n_1,n_2}.
\end{aligned}
\end{equation*}      
It remains to identify the conditions on $\lambda$ under which there hold $\mtx{\Phi}_{12}\geq \mtx{0}$ and $\mtx{\Lambda} \succeq \mtx{0}$ with high probability.

\paragraph{Conditions for $\mtx{\Phi}_{12}> \mtx{0}$:}~Since $\E \mtx{\Phi}_{12} = (\lambda - q) \mtx{J}_{n_1,n_2}$, we see that $\lambda > q$ implies $\E \mtx{\Phi} > \mtx{0}$. In order to guarantee $\mtx{\Phi}_{12} > \mtx{0}$, one may make use of the Chernoff's inequality and the assumption that $n_1 = O(n)$ and $n_2 = O(n)$, which guarantee that  $\mtx{\Phi}_{12} > \mtx{0}$ with high probability as long as $\lambda >q+O\left((\log n)/{n}+\sqrt{(q\log n)/{n}}\right)$.

\paragraph{Conditions for $\mtx{\Lambda} \succeq \mtx{0}$:} Define the projective matrix
\[
\mtx{\Pi}:=\begin{bmatrix}
	\frac{1}{n_1}\mtx{J}_{n_1} & \mtx{0} \\
	\mtx{0} & \frac{1}{n_2}\mtx{J}_{n_2}
\end{bmatrix} 
\]
and let
\[
\widetilde{\mtx{\Lambda}}:= \begin{bmatrix}
		\mbox{diag}(\mtx{A}_{11}\vct{1}_{n_1})-n_1\lambda \mtx{I}_{n_1}-\mtx{A}_{11}+ p \mtx{J}_{n_1}  & -\mtx{A}_{12}+ q \mtx{J}_{n_1,n_2} 
        \\
		- \mtx{A}_{12}^\top + q \mtx{J}_{n_2,n_1} & \mbox{diag}(\mtx{A}_{22}\vct{1}_{n_2})-n_2\lambda \mtx{I}_{n_2}-\mtx{A}_{22} + p \mtx{J}_{n_2} 
	\end{bmatrix} .
\]
The rank-$2$ structure of $\mtx{\Phi}_{12}$ ensures that $\mtx{\Lambda} = (\mtx{I}_n - \mtx{\Pi}) \widetilde{\mtx{\Lambda}}(\mtx{I}_n - \mtx{\Pi})$, where $ \mtx{I}_n $ is the $ n \times n$ identity matrix. Therefore, the condition $\mtx{\Lambda} \succeq \mtx{0}$ is implied by $\widetilde{\mtx{\Lambda}} \succ \mtx{0}$. To proceed, we use the decomposition 
$
\widetilde{\mtx{\Lambda}} = (\E[\mtx{A}] - \mtx{A}) + (\mtx{A} - \E[\mtx{A}] + \widetilde{\mtx{\Lambda}}).
$
It is easy to verify that $\mtx{A} - \E[\mtx{A}] + \widetilde{\mtx{\Lambda}} $ is a diagonal matrix. Letting $n_{\min} := \min (n_1, n_2)$, the Chernoff's inequality ensures that with high probability,
\[
\mtx{A} - \E[\mtx{A}] + \widetilde{\mtx{\Lambda}} \geq (n_{\min}p-2\sqrt{n_{\min}p\log n}-n_{\min}\lambda)\mtx{I}_n.
\]
Moreover, a known result in random matrix theory guarantees that $\|\E[\mtx{A}] - \mtx{A}\| \leq C_0\sqrt{np}$ provided $p \ge {C_0\log n}/{n}$, where $ \| \cdot \| $ denotes the spectral norm; see, e.g., \citet{Feige05}, \citet{Vu2014}, \citet{HajekWuXuSDP14} and \citet{BVH14}. 
We conclude that $\widetilde{\mtx{\Lambda}} \succ \mtx{0}$ as long as $\lambda<p-O(\sqrt{p\log n/ n})$.

The KKT condition in generally does not guarantee $ \mtx{X} $ to be the \emph{unique} solution to \prettyref{eq:SDP1}, establishing which requires additional steps.
We refer the readers to \cite{Cai2014robust} for the details of such steps using Cauchy's interlacing theorem. Summarizing the above arguments, we obtain the following result:
\begin{theorem}
Under the two-cluster $(p,q)$-model with $n_1=O(n)$ and $n_2=O(n)$, if  $p-q > C_0\sqrt{p\log n/ n}$, then there exists a tuning parameter $\lambda$ such that the ground truth partition matrix $ \mtx{X}^\star $ is the unique solution to the convex relaxation \prettyref{eq:SDP1} with high probability.
\end{theorem}


It is worth highlighting that the above argument is quite generic, and can be extended in a straightforward manner to more general and realistic setups of \SBM, such as those with significantly unbalanced clusters, a fast growing number of communities, and general connectivity probability matrix satisfying strong assortativity; see, e.g., \cite{CSX2014}.

%% file: sharp.tex
\section{Sharp threshold for Strong Consistency}
\label{sec:sharp_threshold}
A line of recent work on community detection studies the necessary and sufficient conditions on $p$ and $q$ with \emph{sharp} constants for achieving strong consistency using convex relaxations. This line of work, initiated by \cite{Abbe14} and followed by \citet{HajekWuXuSDP14,Bandeira15,HajekWuXuSDP15,ABBK,perry2015semidefinite}, has achieved remarkable progress. Some assumptions are usually essential to derive conditions with sharp constants, such as relatively balanced clusters and fixed (or slowly growing) number of clusters. Nevertheless, these results demonstrate the strong mathematical power of the convex relaxation approaches, and also highlight the theoretical and practical importance of certain constraints (e.g., row-sum constraints) used in convex relaxations. In this section, we briefly outline the main results in this line, with emphasis on the relevant theoretical tools developed therein. 

As suggested in \citet{AL2014, HajekWuXuSDP15, ABBK}, we here consider an important variation of \prettyref{eq:SDP1} which arises when assuming that all communities have equal size $ n/r $ and that $ r $ is known. In this case, the partition matrix $\mtx{X}$ in \prettyref{eq:membership} satisfies another convex constraint: $\mtx{X} \vct{1}_n = \frac{n}{r} \vct{1}_n$, where $ \vct{1}_n $ is the all one vector in $ \mathbb{R}^n $. Adding this constraint to \prettyref{eq:SDP1} leads to the convex relaxation 
\begin{equation}
\label{eq:SDP_rowsum}
\begin{aligned}
&\text{max} && \left\langle \mtx{X}, \mtx{A} \right\rangle
\\
&\text{subject to} && \mtx{X} \succeq \mtx{0}, 
~ \mtx{X} \geq \mtx{0}, 
~ \mtx{X} \vct{1}_n = \frac{n}{r} \vct{1}_n,  
~ X_{ii} =1 \text{~for~} 1\leq i \leq n.
\end{aligned}
\end{equation}
The term $\lambda \mtx{J}_n$ in \prettyref{eq:SDP1} disappears because $\langle \mtx{X}, \mtx{J}_n \rangle$ is a constant under the row-sum constraint $\mtx{X} \vct{1}_n = \frac{n}{r} \vct{1}_n$. 

We consider the two-cluster $(p,q)$-\SBM~setting as in \prettyref{sec:dual}. For ease of exposition, we further assume that the two clusters are of equal size ($n_1=n_2=n/2$) and that $p = \frac{a\log(n)}{n}$ and $q = \frac{b\log(n)}{n}$, where $a > b > 0$ are two fixed constants. 
In this setting, the information-theoretically limit (i.e., sufficient and necessary condition) for exact recovery is known to be $\sqrt{a}-\sqrt{b}>\sqrt{2}$~\citep{Abbe14,Mossel14}. It turns out that the convex relaxation \prettyref{eq:SDP_rowsum} achieves this information-theoretic limit. In particular, if $\sqrt{a}-\sqrt{b}>\sqrt{2}$, then the convex relaxation \prettyref{eq:SDP_rowsum} with row-sum constraint recovers the true partition matrix $\mtx{X}^*$ in \prettyref{eq:primal_two} as the unique optimum with probability $1-n^{-\Omega(1)}$. Below we present the proof of this result, which makes use of the dual analysis argument.

The dual certificate construction is similar to that introduced in the previous section. A key difference is that the extra row-sum constraint induces a new dual variable $\vec{\mu} \in \reals^n$, allowing for extra freedom in the construction. 
The KKT condition of the relaxation \prettyref{eq:SDP_rowsum} with primal solution $ \mtx{X}^\ast $ defined in \prettyref{eq:primal_two} reads
\begin{itemize}
    \item Stationarity: $-\mtx{A}-\mtx{\Lambda}-\mtx{\Phi}- 
    \mbox{diag}(\vct{\beta})- \vct{\mu} \vct{1}_n^\top - \vct{1}_n
    \vct{\mu}^\top= \mtx{0}$,
	\item Complementary slackness: $\langle \mtx{X}^\star, \mtx{\Lambda}\rangle=\langle \mtx{X}^\star, \mtx{\Phi}\rangle=0$,
    \item Dual feasibility: $\mtx{\Phi} \ge \mtx{0}, \mtx{\Lambda} \succeq \mtx{0}$. 
\end{itemize}
%
%
Under the balanced-cluster assumption $n_1 = n_2 = n/2$, we can simply let $\mtx{\Phi} = \mtx{0}$. Similarly to the argument in \prettyref{sec:dual}, we construct the other dual variables as
\[
\vct{\mu}_1 = - \frac{1}{n_2} \mtx{A}_{12} \vct{1}_{n_2}  + \frac{  \vct{1}_{n_1}^\top \mtx{A}_{12}  \vct{1}_{n_2} }{2 n_1 n_2} \vct{1}_{n_1}, \quad \vct{\mu}_2= - \frac{1}{n_1} \mtx{A}_{12}^\top \vct{1}_{n_1}  +  \frac{\vct{1}_{n_1}^\top \mtx{A}_{12}  \vct{1}_{n_2} }{2 n_1 n_2} \vct{1}_{n_2},
\]
\[
 \vec{\beta_{1} } = - \mtx{A}_{11} \vct{1}_{n_1}  + \mtx{A}_{12} \vct{1}_{n_2}, \quad \vec{\beta_{2} } = - \mtx{A}_{22} \vct{1}_{n_2} + \mtx{A}_{12}^\top \vct{1}_{n_1},
\]
and
\begin{align*}
\mtx{\Lambda} &= -\mtx{A}-\mtx{\Phi}- 
    \mbox{diag}
    (\vct{\beta})- \vct{\mu} \vct{1}_n^\top - \vct{1}_n
    \vct{\mu}^\top  
    \\
    &= (\mtx{I} - \mtx{\Pi})
    \big( (\E[\mtx{A}] - \mtx{A}) + p \mtx{I} - \diag(\vct{\beta}) \big)(\mtx{I} - \mtx{\Pi}),
\end{align*}
where the projection $\mtx{\Pi}$ is defined in \prettyref{sec:dual}. To verify that $\mtx{\Lambda} \succeq \mtx{0}$, it suffices to show that $(\E[\mtx{A}] - \mtx{A}) + p \mtx{I} - \diag(\vct{\beta}) \succ \mtx{0}$. This condition can be established by combining the following two facts. First, by a tight Chernoff's inequality \citep[Lemma 1]{HajekWuXuSDP14}, the inequality $\min_{i\in[n]} (-\beta_i ) \ge \frac{\log n}{\log \log n}$ holds with probability at least $1 - n^{1-(\sqrt{a}-\sqrt{b})^2/2+o(1)}$, which is $1 - n^{-\Omega(1)}$ as long as $\sqrt{a}-\sqrt{b}>\sqrt{2}$. Second, as shown in \prettyref{sec:dual}, $\|  \mtx{A} - \expect{\mtx{A}}\| \leq c' \sqrt{\log n}$ with high probability for a positive constant $c'$. With an addition step proving the uniqueness of the optimal solution, which we skip here, we obtain the following result in the line of \cite{Abbe14,HajekWuXuSDP14,Bandeira15,HajekWuXuSDP15,ABBK}:
\begin{theorem}
Under the two-cluster $(p,q)$-model with equal cluster sizes, suppose that $p = \frac{a\log(n)}{n}$ and $q = \frac{b\log(n)}{n}$. If $\sqrt{a}-\sqrt{b}>\sqrt{2}$, then with  probability $ 1-n^{-\Omega(1)},$ the ground truth partition matrix $ \mtx{X}^\star $ is the unique solution to the convex relaxation~\prettyref{eq:SDP_rowsum}.
\end{theorem}

In the above proof with two clusters, we simply choose $\mtx{\Phi} = \mtx{0}$. This  simple choice not longer works for $r>2$ equal-sized clusters. In that case, denoting by $\mtx{\Phi}_{k,\ell}$ the block in $\mtx{\Phi}$ corresponding to $\mtx{A}_{k, \ell}$, we may choose $\mtx{\Phi}_{k,\ell}=\vct{y}_{k,\ell} \vct{1}^\top_{n_\ell} +  \vct{1}_{n_k} \vct{z}^\top_{k,\ell}$ for some vectors $\vct{y}_{k,\ell}$ and $\vct{z}_{k,\ell}$. The detailed argument can be found in \cite{HajekWuXuSDP15}.

We note that the dual variable $-\beta_i$ above corresponds to the number of neighbors  node $i$ has in its own cluster minus the number of its neighbors in the other cluster. The variable $-\beta_i$ is closely related to the information-theoretic lower bounds of strong consistency. In particular, it was shown in~\cite{Abbe14}  that if $\sqrt{a}-\sqrt{b} < \sqrt{2}$,  then with probability $1-o(1)$, there exists a pair of nodes $i$ and $j$ from different clusters such that $-\beta_i < -1$ and $-\beta_j < -1$, which further implies that the maximum likelihood estimator cannot coincide with the ground truth $\mtx{X}^\star$ and hence the impossibility of achieving strong consistency.

%
%
%

\textbf{Other convex relaxation formulations:}  It is noteworthy that besides the formulation~\prettyref{eq:SDP_rowsum} with row-wise constraints, there exist other convex optimization methods that achieve the sharp thresholds for strong consistency. For example, in the two-community setting, one may consider the following formulation:
\begin{equation}
\label{eq:SDP_centered}
\begin{aligned}
&\text{max} && \left\langle \mtx{Y}, \mtx{A}-\lambda \mtx{J}_n \right\rangle
\\
&\text{subject to} && \mtx{Y} \succeq \mtx{0}, ~Y_{ii} =1 \text{~for~} 1\leq i \leq n,
\end{aligned}
\end{equation}
which can be seen as a variant of the well-known SDP relaxation for MaxCut \citep{goemans1995improved}.
With some appropriate $\lambda$, the formulation~\prettyref{eq:SDP_centered} is shown to recover the ground truth centered partition matrix $\mtx{Y}^* = 2\mtx{X}^* - \mtx{J}_n$ as the unique optimum with probability $1-o(1)$, within the information-theoretically feasible range of \SBM~parameters; see \citet{Abbe14, HajekWuXuSDP14,Bandeira15} for the case of equal cluster sizes, and \citet{HajekWuXuSDP15} for the unbalanced case. Extensions to the more general setting with multiple clusters of unequal sizes can be found in \cite{perry2015semidefinite}.

%% file: weak_assortativity.tex
\section{Projective Matrix Based Convex Optimization and Weak Assortativity}
\label{sec:weak_assortativity}
In \prettyref{sec:dual}, we present strong consistency results for the relaxation \prettyref{eq:SDP1} under the $(p, q)$-\SBM, which can be extended to more general connectivity probability matrices $\mtx{B}$. However, the strong consistency of \prettyref{eq:SDP1} requires that the tuning parameter $\lambda$ is between minimum within-community edge density $\min_{1 \leq a \leq r} B_{aa}$ and the maximum cross-community edge density $\max_{1 \leq  a < b \leq r} B_{ab}$. 
This means that one must have  $\min_{1 \leq a \leq r} B_{aa}> \max_{1 \leq  a < b \leq r} B_{ab}$, namely, the strong assortativity of $\mtx{B}$. The notion of strong assortativity, coined in \cite{AL2014}, is a common assumption in the literature of community detection under the \SBM; see, e.g., \cite{RCY2011} and \cite{CCT2012}.

Strong assortativity is sometimes too restrictive. To address this issue, \cite{AL2014} consider the \SBM under the weak assortativity assumption $B_{aa} > \max_{b \neq a} B_{ab}$ for any $a = 1, \ldots, r$, and show that the convex relaxation~\prettyref{eq:SDP_rowsum} achieves strong consistency under certain conditions of the model parameter. 
The dual analysis used to prove this result is similar to that in \prettyref{sec:sharp_threshold}, although they do not focus on identifying sharp thresholds. However, this result relies on the strong assumption that all cluster have equal sizes. Indeed, as discussed in \prettyref{sec:sharp_threshold}, this assumption is necessary for the true partition matrix $ \mtx{X}^\ast $ to be feasible to the formulation~\prettyref{eq:SDP_rowsum}.

It is desirable to develop a convex relaxation that does not require this unrealistic  assumption of equal-cluster-size while still guaranteeing strong consistency under weak assortativity. To this end, we keep a version of the row-sum constraint $  \mtx{X} \vct{1}_n = \frac{n}{r} \vct{1}_n $, which appears essential for the weak assortativity setting, but try to remove its explicit nal constraints in \prettyref{eq:SDP_rowsum} with a trace constraint, i.e.,
\begin{equation*}
\begin{aligned}
&\max_{\mtx{X}} && \left\langle \mtx{X}, \mtx{A} \right\rangle
\\
&\text{subject to} && \mtx{X} \succeq \mtx{0}, ~ \mtx{X} \geq \mtx{0}, ~ \mtx{X} \vct{1}_n = \frac{n}{r} \vct{1}_n, ~\text{trace}(\mtx{X})=n.
\end{aligned}
\end{equation*}
With a change of variable $\mtx{Z} = \frac{r}{n}\mtx{X}$, we obtain the following relaxation:
\begin{equation}
\label{eq:SDP_PW}
\begin{aligned}
&\max_{\mtx{Z}} && \left\langle \mtx{Z}, \mtx{A} \right\rangle
\\
&\text{subject to} && \mtx{Z} \succeq \mtx{0},~\mtx{Z} \geq \mtx{0},~ \mtx{Z} \vct{1}_n = \vct{1}_n,~\text{trace}(\mtx{Z}) = r.
\end{aligned}
\end{equation}
Note that the above relaxation has no explicit dependence on the cluster sizes. This relaxation appears less studied in the community detection literature compared to \prettyref{eq:SDP1} and \prettyref{eq:SDP_rowsum}, but it in fact coincides with the well-known Peng-Wei relaxation \citep{PengW07} for K-means in the literature of clustering in Euclidean space. The work in that literature shows that under certain affinity or separation conditions, the ground truth projective matrix
\begin{equation}
\label{eq:projection}
\mtx{Z}^\star=\mtx{P}
\begin{bmatrix} 
\frac{1}{n_1}\mtx{J}_{n_1} &~ & ~
\\ 
~ &\ddots & ~
\\ 
~ & ~ & \frac{1}{n_r}\mtx{J}_{n_r}
\end{bmatrix}
\mtx{P}^{\top},
\end{equation}
instead of the true partition matrix $ \mtx{X}^\star $ in the form of \prettyref{eq:membership}, is the unique solution to the relaxation \prettyref{eq:SDP_PW}; see, e.g., \cite{iguchi2015tightness, IguchiMPV15}, \cite{LLLSW2017} and \cite{fei2018hidden}. 

The relaxation \prettyref{eq:SDP_PW} was recently analyzed in the context of community detection by \cite{YSC2018}. They show that this relaxation enjoys strong consistency (in the sense of recovering $ \mtx{Z}^\star $ as the unique solution with high probability) under \SBM~with weak assortativity and unbalanced communities. This result substantially extends the result in \cite{AL2014}.

The proof of the above result again uses dual analysis. The primal feasibility of the projective matrix $ \mtx{Z}^\star $ in \prettyref{eq:projection} to the relaxation \prettyref{eq:SDP_PW} can be directly verified. To further establish the optimality of $ \mtx{Z}^\star $ to \prettyref{eq:SDP_PW}, it suffices, similarly to \prettyref{sec:dual} and \prettyref{sec:sharp_threshold}, to construct dual variables $\mtx{\Phi}$, $\mtx{\Lambda}$, $\vct{\mu}$ and $z$ for which the following KKT conditions hold:
\begin{itemize}
    \item Stationarity: $-\mtx{A}-\mtx{\Lambda}-\mtx{\Phi}+ z \mtx{I} - \vct{\mu} \vct{1}_n^\top - \vct{1}_n
    \vct{\mu}^\top= \mtx{0}$,
	\item Complementary slackness: $\langle \mtx{Z}^\star,  \mtx{\Lambda}\rangle=\langle\mtx{Z}^\star , \mtx{\Phi}\rangle=0$,
    \item Dual feasibility: $\mtx{\Phi} \ge \mtx{0}, \mtx{\Lambda} \succeq \mtx{0}$. 
\end{itemize}
Notice that the above conditions differ from the KKT conditions of the relaxation \prettyref{eq:SDP_rowsum} given in \prettyref{sec:sharp_threshold}; in particular, the diagonal matrix $-\diag(\vct{\beta})$ is replaced by the scaled identity matrix $z\mtx{I}$, which reduces the freedom for the dual construction. The detailed construction of these dual variables can be found in \cite{YSC2018}. Here we cite the result established therein:

\begin{theorem}
Denote by $n_{\max}$ and $n_{\min}$ the maximum and minimum cluster sizes in the~\SBM, respectively. If the following separation condition holds:
\begin{align*}
\min_k \left(B_{kk} - \max_{l \neq k} B_{kl}\right) \geq& 2\sqrt{6\log n} \max_k \sqrt{B_{kk}/n_k} + 6 \max_{1 \leq k < l \leq r} \sqrt{B_{kl} \log n / n_{min}} 
\\
&+ C_0\sqrt{(n/n_{\min}^2)(\max_k B_{kk})},
\end{align*}
then with high probability, the matrix $ \mtx{Z}^\star  $ in \prettyref{eq:projection} is the unique solution to \prettyref{eq:SDP_PW}.
\end{theorem}

%% file: primal.tex
\newcommand{\Xhat}{\ensuremath{ \widehat{\mtx{X}} }}
\newcommand{\Xstar}{\ensuremath{ \mtx{X}^\star }}
\newcommand{\xhat}{\ensuremath{ \widehat{X} }}
\newcommand{\xstar}{\ensuremath{ X^* }}
\newcommand{\onenorm}[1]{\ensuremath{ \left\| #1 \right\|_{1} }}
\newcommand{\error}{\ensuremath{\gamma}}
\newcommand{\defn}{:=}
\newcommand{\onevec}{\vct{1}}
\newcommand{\OneMat}{\mtx{J}_n}
\newcommand{\size}{\ell}
\newcommand{\numclust}{r}

\section{Weak Consistency by Primal Analysis}
\label{sec:primal}

In this section, we describe a ``primal'' approach for analyzing convex relaxation formulations. As opposed to the ``dual'' analysis presented in the previous sections, this approach directly makes use of the primal optimality and feasibility of the SDP solution $\widehat{\mtx{X}}$ and the ground-truth $\Xstar$. Combined with the celebrated Grothendieck's inequality~\citep{Grothendieck}, this approach leads to a simple proof of weak consistency~\citep{Vershynin14}. Under the equal-cluster-size assumption, a more refined primal analysis provides unified guarantees covering both weak and strong consistency~\citep{FeiChen2017_exponential}.

To present the above results, we consider the $ (p,q) $-\SBM and the convex relaxation~\eqref{eq:SDP1} with $q < \lambda < p$. Let $\tau := \min(p-\lambda, \lambda - q)$. For expository convenience, we simply choose $ \lambda = \frac{p+q}{2}$, which is an approximation of~\eqref{eq:lambda_p_q} used in the exact log likelihood. In this case we have $\tau = \frac{p - q}{2}$. The arguments below can be extended to other values of $ \lambda $ in a straightforward fashion.
Recall that our goal is to estimate the ground truth partition matrix  $\mtx{X}^\star$. We use the solution $\widehat{\mtx{X}}$ to the convex relaxation \eqref{eq:SDP1} as an estimator, and aim to characterize the efficiency of \eqref{eq:SDP1} by bounding the distance between $\widehat{\mtx{X}}$ and $\Xstar$.

The primal analysis begins with noting that the optimality of $\widehat{\mtx{X}}$ and the feasibility of $\Xstar$ to the convex relaxation \eqref{eq:SDP1} imply that $\left\langle \mtx{X}^*, \mtx{A}-\lambda \mtx{J}_n \right\rangle \leq \left\langle \widehat{\mtx{X}}, \mtx{A}-\lambda \mtx{J}_n\right\rangle$. 
Rearranging the inequality and separating the expectation and deviation terms, we obtain the following inequality:
\begin{align}
0 & \le\left\langle \Xhat-\Xstar,\mtx{A}-\lambda \OneMat \right\rangle \nonumber  =\left\langle \Xhat-\Xstar,\E\mtx{A}-\lambda \OneMat \right\rangle +\left\langle \Xhat-\Xstar,\mtx{A}-\E\mtx{A}\right\rangle .\label{eq:basic_inequality}
\end{align}
To proceed, we make use of a simple observation on the relationship between $\mtx{X}^\star$ and $\E \mtx{A}$. For each $i \neq j$, $\E A_{ij} =p  > (p+q)/2$ if and only if $X_{ij}^*=1$, whereas $\E A_{ij} = q < (p+q)/2$ if and only if $X_{ij}^*=0$. Due to the feasibility of $\widehat{\mtx{X}}$, we deduce that the entries of the error matrix $\Xstar-\Xhat$ must have matching signs with those of $\E\mtx{A}-\frac{p + q }{2}\OneMat$. This observation implies that
\begin{equation*}
	\left\langle \Xstar-\Xhat,\E\mtx{A}-\lambda \OneMat \right\rangle \ge\frac{p - q }{2} \|\Xhat - \Xstar\|_1.\label{eq:expectation_gamma}
\end{equation*}


The key step in~\cite{Vershynin14} involves bounding the ``deviation term'' $ \langle \Xhat-\Xstar,\mtx{A}-\E\mtx{A} \rangle  $. To this end, one applies the triangle inequality to obtain
\begin{align*}
\left\langle \Xhat-\Xstar,\mtx{A}-\E\mtx{A}\right\rangle 
&\leq |\langle \Xhat, \mtx{A}-\E\mtx{A} \rangle| + \left|\left\langle \Xstar, \mtx{A}-\E\mtx{A} \right\rangle \right| \leq 2\sup_{\substack{\mtx{X}\succeq0 \\ \textmd{diag}(\mtx{X})\le\onevec}} \left|\left\langle \mtx{X}, \mtx{A} - \E\mtx{A}\right\rangle \right|,
\end{align*}
where $\textmd{diag}(\mtx{X})\le\onevec$ means that $X_{ii} \leq 1$ for all $i=1, \ldots, n$, and the second inequality follows from the feasibility of $\Xhat$ and $\Xstar$. The Grothendieck's inequality~\citep{Grothendieck,Lindenstrauss}
guarantees that 
\[
\sup_{\mtx{X}\succeq0,\textmd{diag}(\mtx{X})\le\onevec}\left|\left\langle \mtx{X},\mtx{A}-\E\mtx{A}\right\rangle \right|
\leq K_{G}\norm{\mtx{A}-\E\mtx{A}}_{\infty\to1}
\]
where $K_{G}$ denotes the Grothendieck's constant ($0<K_{G}\leq1.783$)
and $\norm{\bM}_{\infty\to1} \defn \sup_{\vct{x}:\|\vct{x}\|_{\infty}\leq1}\norm{\bM\vct{x}}_{1}$
is the $\ell_{\infty}\to\ell_{1}$ operator norm for a matrix $\bM$.
We next make use of the identity 
\begin{align*}
\norm{\mtx{A}-\E\mtx{A}}_{\infty\to1} &=\sup_{\substack{\vct{x}:\|\vct{x}\|_{\infty}\leq 1 \\ \vct{y}:\|\vct{y}\|_{\infty}\leq 1}} \vct{y}^{\top}(\mtx{A}-\E\mtx{A})\vct{x} =\sup_{\substack{\vct{x}\in\{\pm1\}^{ n } \\ \vct{y} \in\{\pm1\}^{ n } }} \vct{y}^{\top}(\mtx{A}-\E\mtx{A})\vct{x}.
\end{align*}
For each pair of fixed vectors $\vct{x},~\vct{y}\in\{\pm1\}^{n}$, $ \vct{y}^{\top}(\mtx{A}-\E\mtx{A})\vct{x}$ can be written as the sum of independent random variables and hence be bounded using Bernstein's inequality. Taking a union bound over all $ \vct{x} $ and $ \vct{y} $, we obtain that with probability at least $1- 2(e/2)^{-2 n }$,
\[
\norm{\mtx{A}-\E\mtx{A}}_{\infty\to1} \leq 2\sqrt{2 p ( n ^{3}- n ^{2})}+ ({8}/{3}) n .
\]

Combining pieces, we conclude that with probability at least $1-2(e/2)^{-2 n }$,
\begin{align*}
\|\Xhat-\Xstar\|_1 & \leq {2}/{( p - q) }\left(8\sqrt{2 p ( n ^{3}- n ^{2})}+ ({32}/{3}) n \right)\overset{(i)}{\le}{45\sqrt{ p  n ^{3}}}/{( p - q) },
\end{align*}
where step $(i)$ holds provided that $ p \geq1/ n $. We have therefore established the following result, which was first proved in \cite{Vershynin14}:
\begin{theorem}
\label{thm:primal}
Under the $(p,q)$-SBM, if $ p \ge 1/n $ and $p>q$, then with high  probability, any optimal solution $ \Xhat $ of the SDP~\eqref{eq:SDP1} satisfies the bound $\|\Xhat - \Xstar\|_1 \le {45\sqrt{ p  n ^{3}}}/{ (p - q) }$.
\end{theorem}

A prominent advantage of this result is that it is applicable even when the networks is very sparse, i.e., $ p \asymp 1/n $. Moreover, the result does not require any assumption on the number of clusters, the cluster sizes, or the knowledge thereof. 

Given the above bound on $\|\Xhat - \Xstar\|_1$, one can further derive bounds on the clustering errors. For example, in the two-cluster case one can simply estimate the underlying community label $\phi$ by taking the entry-wise sign of the leading eigenvector of $\widehat{\mtx{Y}} = 2\Xhat - \mtx{J}_n$, the centered version of $ \Xhat $.
It is easy to show that the fraction of misclassified vertices is upper bounded by $(1/n^2)\|\Xhat - \Xstar\|_1$ up to constant factors~\citep{Vershynin14}. Combined with \prettyref{thm:primal}, we conclude that the fraction of misclassified vertices is smaller than $\epsilon$ with high probability provided that $n (p-q)^2 \gtrsim p \epsilon^{-2} $.

Under the equal-cluster-size assumption, the error bound in \prettyref{thm:primal} has been substantially improved by~\cite{FeiChen2017_exponential}. In particular, by deriving a tighter bound on the term $ \langle \Xhat-\Xstar,\mtx{A}-\E\mtx{A} \rangle  $, they show that the estimation error of \prettyref{eq:SDP1} in fact decays exponentially:
\begin{align}
\label{eq:exponential}
\|\Xhat - \Xstar\|_1 \lesssim n^2 \exp\left[-\Omega\left( \frac{n(p-q)^{2}}{p\numclust} \right)\right]
\end{align}
provided that $ {(p-q)^{2}}/{p} \gtrsim {\numclust^2}/{n}. $ Assuming $r=O(1)$, the bound~\eqref{eq:exponential} guarantees that the fraction of misclassified nodes is at most $ \epsilon $ as long as $n(p-q)^{2} \gtrsim p \log (\epsilon^{-1}) $. Consequently, the convex relaxation achieves weak consistency (i.e., $ \epsilon \to 0 $) if $n{(p-q)^{2}}/{p\log n} \to \infty  $, and strong consistency (i.e., $ \epsilon <1/n $, which implies $ \epsilon =0 $) if  $n{(p-q)^{2}}/{p} \gtrsim \log n  $.

Let us provide an additional remark on the setting with two equal-sized clusters and sparse networks, i.e., $p=a/n$ and $q=b/n$ for two constants $a,b$. Assuming $a>b$ and $ (a-b)^2 \ge C_0 (a+b) $, a related question is what is the smallest $C_0$ to guarantee that the SDP formulation can achieve $ \epsilon<1/2 $, that is,  producing a non-trivial estimate that is better than random guess. It is known that the information-theoretic \emph{necessary} condition for achieving $ \epsilon < 1/2 $ fraction of misclassified nodes is $(a-b)^2> 2(a+b)$ \citep{Mossel12,Massoulie13,Mossel13}. In the converse direction, recent work by \citet{MontanariSen15} shows that the SDP \eqref{eq:SDP_centered} with appropriate $\lambda$ achieves $ \epsilon<1/2 $ with probability $1-o_n(1)$ if $(a-b)^2 > (2+\epsilon) (a+b)$ and $a+b >C_1$ for a sufficiently large constant $ C_1 $ depending on $\epsilon$. Therefore, the results in \cite{MontanariSen15} show that the convex relaxation approach is nearly optimal in achieving non-trivial estimation. More recent work by~\cite{FanMontanari16} gives more precise bounds on $C_1$ by proving that a simple local algorithm approximately solves the convex program within a factor $1+O(1/(a+b)).$

%% file: cmm.tex
\section{Degree Correction and Convexified Modularity Maximization}

 An alternative to the likelihood-based approach~\prettyref{eq:likelihood} for community detection is modularity maximization. Proposed by Newman and Girvan \citep{Newman2006modularity}, this approach involves finding a partition matrix $ \mtx{X} $ that maximizes the so-called ``modularity'':
 \begin{equation}
\label{eq:modularity}
\max_{\mtx{X}} \sum_{1 \leq i < j\leq n}\left(A_{ij}-\frac{d_{i}d_{j}}{2L}\right)X_{ij}
\end{equation}
where $d_{i}:=\sum_{j=1}^n A_{ij}$ is the degree of node $i$, and $L := \frac{1}{2}\sum_{i=1}^n d_{i}$ is the total number of edges. Compared to likelihood maximization~\prettyref{eq:likelihood}, modularity maximization can be viewed as an adaptive variant, where the single tuning parameter $\lambda$ is replaced by the quantity $\frac{d_{i}d_{j}}{2L}$ that is adaptive to each pair of nodes $(i, j)$.

To overcome the so-called ``resolution limit'' of modularity maximization \citep{fortunato2007resolution}, the work in \citet{RB2006,LF2011} proposed replacing the number $\frac{1}{2L}$ in \prettyref{eq:modularity} with a tuning parameter $\lambda$. Convexifying the program similarly to \prettyref{eq:SDP1}, we obtain the following convex relaxation:
\begin{equation}
\label{eq:cvx}
\begin{aligned}
&\max_{\mtx{X}} && \iprod{ \mtx{X} }{ \mtx{A}- \lambda \vct{d}\vct{d}^\top }
\\
&\text{subject to} && \mtx{X} \succeq \mtx{0},~\mtx{X} \geq \mtx{0},~X_{ii} =1, \text{~for~} 1\leq i \leq N.
\end{aligned}
\end{equation}
This approach, termed \emph{convexified modularity maximization}, was first proposed by~\citet{CLX2015}. Notice that if we replace $\vct{d}$ with $\vct{1}_n$, then the formulation~\prettyref{eq:cvx} falls back to convexified maximum likelihood~\prettyref{eq:SDP1}.

An important advantage of modularity maximization over likelihood maximization lies inits applicability beyond the standard \SBM. Indeed, it was shown in~\citet{CLX2015} that the formulation~\prettyref{eq:cvx} is consistent even if the data is generated from the more general \emph{degree corrected} stochastic block model (\DCSBM). We now describe this result, beginning with an introduction of \DCSBM.

One obvious limitation of standard \SBM~is that nodes within the same group are statistically equivalent and hence have homogeneous degrees. Therefore, standard \SBM~fails to model degree heterogeneity that is common in real world networks. This limitation motivates researchers to consider the more general \DCSBM. In \DCSBM, each pair of nodes $i$ and $j$ are connected with probability $\theta_i \theta_j B_{\phi(i)\phi(j)}$---as opposed to $B_{\phi(i)\phi(j)}$ in \SBM---where $\theta_i$ is a degree heterogeneity parameter for the node $i$.

A first step in deriving statistical guarantees under \DCSBM~is an appropriate generalization of the ``density gap'' condition $p-q>0$ used in the $(p, q)$-\SBM. We first note that for each node $i$ in community $a$, its expected degree is
\begin{align*}
\E d_i & = \sum_{j \neq i} \E A_{ij} = \sum_{j \neq i} \theta_i \theta_j B_{a \phi(j)} = \theta_i \sum_{j} \theta_j B_{a \phi(j)} - \theta_i^2 B_{aa} \approx \theta_i \sum_{j} \theta_j B_{a \phi(j)}
\end{align*}
Define $H_a := \sum_{j} \theta_j B_{a \phi(j)}$ and note the alternative expression $H_a = \sum_{b=1}^r B_{ab}G_b$, where $G_a := \sum_{i : \phi(i) = a} \theta_i $. Here $H_a$ captures the \emph{average degree} of the nodes in the $a$-th cluster in the sense that $H_a \approx \E d_i / \theta_i$. Based on these observations, \citet{CLX2015} proposed the following \emph{degree-corrected density gap condition}:
\begin{equation}
\label{eq:dcgap}
\max_{1 \leq a < b \leq r} \frac{B_{ab}}{H_a H_b} < \min_{1 \leq a \leq r} \frac{B_{aa}}{H_a^2}.
\end{equation}
This degree-corrected density gap condition enjoys several desirable properties. First, it only depends on the community aggregated quantities $H_a$'s instead of individual node parameters $\theta_i$'s, so it is robust against abnormal nodes, such as one small $\theta_i$. Moreover, while different pairs of $(\mtx{B}, \vct{\theta})$ may correspond to the same \DCSBM, the condition \prettyref{eq:dcgap} is invariant under equivalent \DCSBM's.

\citet{CLX2015} showed that under appropriate conditions stated in terms of the gap condition~\prettyref{eq:dcgap}, the convexified modularity maximization \prettyref{eq:cvx} provides a good estimate of the ground truth partition matrix $\mtx{X}^*$:

\begin{theorem}
\label{thm:approx}
Under \DCSBM, suppose that the degree-corrected density gap condition \prettyref{eq:dcgap} holds and the tuning parameter $\lambda$ in the formulation~\prettyref{eq:cvx} satisfies
\begin{equation}
\label{eq:RDGC}
\max_{1 \leq a < b \leq r} \frac{B_{ab} + \delta}{H_a H_b} \le \lambda \le \min_{1 \leq a \leq r} \frac{B_{aa} - \delta}{H_a^2}
\end{equation}
for some number $\delta>0$. Then with high probability, any solution~$ \widehat{ \mtx{X} } $ to the convex relaxation~\prettyref{eq:cvx} satisfies the bound
\begin{equation}
\label{eq:pillar1}
\sum_{i=1}^{n}\sum_{j=1}^{n} \theta_i \theta_j |\widehat{X}_{ij} - X_{ij}^*| \leq \frac{C_0}{\delta}\left( 1+ \lambda  \sigma \right)\left(  \sqrt{n \sigma }  + n \right),
\end{equation}
where  $\sigma := \sum_{a, b} B_{ab} G_a G_b$.
\end{theorem}

This theorem is proved via a primal analysis argument similar to that in~\citet{Vershynin14} as presented in \prettyref{sec:primal}. 

Note that the error metric in \prettyref{eq:pillar1} is weighted by the degree heterogeneity parameters. This captures the fact that nodes with different $ \theta_i $'s and hence different degrees have different contributions to the overall clustering quality. Moreover, the error bound~\prettyref{eq:pillar1} is insensitive to $\theta_{\min} :=\min_i \theta_{i}$. Therefore, the presence of a small $ \theta_{\min} $ does not hinder recovery of the memberships of nodes with large $\theta_i$. 

To better appreciate the results in \prettyref{thm:approx},
let us consider a simple sub-class of \DCSBM with symmetric and balanced clusters, where $B_{aa} = p$ for all $1 \le a \le r$, $B_{ab}=q$ for all $1\leq a<b \leq r$, and $G_a = g$ for all $ 1\le a \le r $.
In this setting, direct calculation gives $H_a = ((r-1)q + p)g$ for all $1 \le a \le r$. The constraint \prettyref{eq:RDGC} on $\delta$ and $\lambda$ then simplifies to
\begin{equation}
\label{eq:RDGC2}
p-q \ge 2\delta \quad \text{and} \quad \frac{q+\delta}{(p+(r-1)q)^2g^2} \le \lambda  \le \frac{p - \delta}{(p+(r-1)q)^2g^2}.
\end{equation}
Note that the first inequality above is identical to the standard density gap condition imposed in \SBM. Furthremore, one has $\sigma=r(p+(r-1)q)g^2 \le r^2 p g^2$. Substituting these expressions into equation~\prettyref{eq:pillar1}
yields the error bound
\begin{equation}
\label{eq:pillar2}
\sum_{i=1}^{n}\sum_{j=1}^{n} \theta_i \theta_j |\widehat{X}_{ij} - X_{ij}^*|  \lesssim  \frac{r}{\delta}  \left( rg\sqrt{np} + n \right).
\end{equation}
As discussed in the previous sections, 
the solution $ \Xhat $ of the convex relaxation is not necessarily a partition matrix;
\cite{CLX2015} show that an explicit clustering can be extracted from $ \Xhat $ using a \emph{weighted} $k$-medoids algorithm. The number of misclassified nodes, \emph{weighted} by their degree heterogeneity parameters $ \{ \theta_i \}$, satisfies the bound
\begin{align}
\label{eq:Sbound}
\sum_{i \in \mathcal{E} } \theta_{i} \lesssim \frac{r}{\delta} \left(  r\sqrt{np} + \frac{n}{g} \right),
\end{align}
where $\mathcal{E}$ is the set of misclassified nodes.
We refer the readers to \citet[Sections 2.4, 3.2]{CLX2015} for the details of this result.

We note that the bound in \prettyref{eq:Sbound} gives non-trivial guarantees all the way down to the \emph{sparse graph regime} with bounded average degrees. For example, suppose that $p =a/n$ and $q=b/n$ for two fixed constants $a>b$,  $r=O(1)$, and $g \asymp n$. With $(a-b)/\sqrt{a}$ sufficiently large and the choice $ \delta \asymp (a-b)/n$, one may bound the right hand side of~\prettyref{eq:Sbound} by an arbitrarily small constant times $ n $, thereby guaranteeing that an arbitrarily small fraction of nodes are misclassified.
In contrast, standard spectral clustering methods are known to be inconsistent in this sparse regime \citep{KMMNSZZ2013}. 

\paragraph{Strong Consistency:}

If in addition the minimum degree heterogeneity parameter $\theta_{\min}$ is not too small, \cite{CLX2015} 
further show that convexified maximum likelihood perfectly recovers the ground truth communities with high probability.
To illustrate this strong consistency result, consider a sub-class of \DCSBM~with $B_{aa} = p$ for all $a=1, \ldots, r$, and $B_{ab}=q$ for all $1\leq a<b \leq r$. Under this setup, the degree-corrected density gap condition \prettyref{eq:RDGC} becomes
\begin{equation}
\label{eq:sep_perf}
\max_{1 \leq a < b \leq r} \frac{q + \delta}{H_a H_b}  \le  \lambda  \le \min_{1 \leq a \leq r} \frac{p - \delta}{H_a^2}.
\end{equation}
Further define $G_{\min}:= \min_{1 \leq a \leq r} G_a.$
The following theorem, again extracted from~\citet{CLX2015}, provides sufficient conditions for strong consistency.
\begin{theorem}
\label{thm:exact}
Suppose that the tuning parameter $\lambda $ satisfies the degree-corrected density gap condition \prettyref{eq:sep_perf} for some number $\delta >0$, and that
\begin{align}
 \label{eq:assump}
 \delta > C_0 \left( \frac{\sqrt{q n}}{G_{\min} }  + \sqrt{\frac{p \log n}{G_{\min}\theta_{\min}}} \right).
\end{align}
With high probability, the ground truth partition matrix $ \mtx{X}^\star $ is the unique solution to the convex relaxation~\prettyref{eq:cvx}.
\end{theorem}

This theorem is proved using a dual analysis argument similar to those presented in \prettyref{sec:dual}. Note that the condition~\prettyref{eq:assump} for perfect recovery depends on the minimum values $G_{\min}$ and $\theta_{\min}$. This is expected: it is impossible to recover the membership of a node with an overly small $ \theta_i $ and $ G_a $, as this node will have too few edges. 

%% file: robustness.tex
\section{Robustness against Outlier nodes}
\label{sec:robustness}

The results presented in the previous sections highlight some distinguishing statistical advantages of convex relaxation methods for community detection, including consistency under sparse networks and weak assortativity. In this section, we demonstrate another important benefit of convex relaxation, namely, robustness. In \citet{Cai2014robust}, it is shown, both theoretically and empirically, that a simple variant of \prettyref{eq:SDP1} is robust against outlier/adversarial nodes. That is, convex optimization methods can yield consistent estimates even if a proportion of the nodes have arbitrary connections that do not satisfy the \SBM. Below we give an expository introduction of this result and the associated analytical arguments.

For analytical convenience,  \citet{Cai2014robust} consider the following modified version of the convex relaxation \prettyref{eq:SDP1}:
\begin{equation}
\label{eq:SDP2}
\begin{aligned}
&\max_{\mtx{X}} && \left\langle \mtx{X}, \mtx{A} - \lambda \mtx{J} - (\alpha - \lambda) \mtx{I} \right\rangle
\\
&\text{subject to} && \mtx{X} \succeq \mtx{0},~\mtx{X} \geq \mtx{0},~ {X}_{ii}\leq 1,\text{~for~} 1\leq i \leq n.
\end{aligned}
\end{equation}
The optimization problem \prettyref{eq:SDP2} differs from \prettyref{eq:SDP1} in that it does not require all diagonal entries of $\mtx{X}$ to equal one. Instead, the entries of $\mtx{X}$ are constrained to fall between $0$ and $1$, and the trace of $\mtx{X}$ is penalized by introducing a second tuning parameter $\alpha$. 

To study the robustness of the formulation~\prettyref{eq:SDP2}, let us consider an extension of \SBM~that allows for abnormal nodes. We assume that the set of nodes $V$ can be decomposed as $V = I \cup O$, where $I$ represents the inliers and $O$ is the outliers. We assume the edges between nodes within $I$ are generated randomly according to the $ (p,q) $-\SBM, whereas if one of $i$ and $j$ is an outlier, then $i$ and $j$ are connected in an arbitrary or even adversarial manner.

Under the above model, one can show that the formulation~\prettyref{eq:SDP2} can still consistently recover the community memberships of the nodes in $I$. The argument is based on analyzing the KKT condition for \prettyref{eq:SDP2}, but it differs from the purely dual analysis presented in Section~\ref{sec:dual}. There when studying the standard \SBM~and the formulations \prettyref{eq:SDP1}, \prettyref{eq:SDP_rowsum} and \prettyref{eq:SDP_PW}, it is clear which is the desired solution whose optimality we want to certify---it is either the ground truth community matrix in \prettyref{eq:membership} or the ground truth projective matrix in \prettyref{eq:projection}. Under the generalized \SBM
with outlier nodes, however, the desired solution to \prettyref{eq:SDP2} is much more complicated. Consequently, the dual analysis of \prettyref{eq:SDP2} involves both the construction of a candidate primal solution and its dual certification. This primal-dual analysis is powerful but relatively uncommon in the convex community detection literature. Below we briefly outline the key ideas of this argument.

Similarly to \prettyref{sec:dual}, we focus on the $(p, q)$-model with two inlier clusters.  We assume that the $ |I| =n $ inliers are partitioned into two clusters of sizes $ n_1 $ and $ n_2 $, and the number of outliers is $|O| = m$, with $n=n_1 + n_2 + m$ being the total number of nodes. We refer to this model as the \emph{$(p, q, n_1, n_2, m)$-generalized~\SBM}. In the presence of outliers, the adjacency matrix (up to some permutation)  has the following augmented block structure:
\begin{equation}
\label{eq:adjacency}
\mtx{A} = \begin{bmatrix} 
\mtx{A}_{11} & \mtx{A}_{12} & \mtx{Z}_1
\\
\mtx{A}_{12}^{\top} & \mtx{A}_{22}    &  \mtx{Z}_2 
\\ 
\mtx{Z}_1^{\top} & \mtx{Z}_2^{\top} & \mtx{W}
\end{bmatrix}.
\end{equation}
Here $\mtx{A}_{11}$, $\mtx{A}_{22}$ and $\mtx{A}_{12}$ correspond to the connectivity within and between the two inlier clusters, and hence have the same statistical properties as in \SBM. 
On the other hand, the blocks $\mtx{Z}_{1}$ and $\mtx{Z}_2$ represent the connectivity between the outliers and inliers, whereas the $ m\times m$ symmetric block $\mtx{W}$ is the adjacency matrix within the outliers. Note that $\mtx{Z}_{1}$, $\mtx{Z}_2$ and $\mtx{W}$ are all arbitrary.
 
We say that the formulation~\prettyref{eq:SDP2} is robust against outliers if it achieves strong consistency within the inlier nodes in $I$; that is, its solution $\widehat{\mtx{X}}$ is of the form
\begin{equation}
\label{eq:output}
\widehat{\mtx{X}}=
\begin{bmatrix} 
\mtx{J}_{n_1} & \mtx{0} & \star
\\ 
\mtx{0} &  \mtx{J}_{n_2}    &  \star 
\\ 
\star & \star & \star 
\end{bmatrix}.
\end{equation}
To establish this property, we begin with the KKT condition of \prettyref{eq:SDP2} for the primal solution $\widehat{\mtx{X}}$:
\begin{equation}
\label{eq:KKT_robust}
\begin{cases}
\text{Stationarity:~}-\mtx{A} + \lambda \mtx{J} + (\alpha - \lambda) \mtx{I} = \mtx{\Lambda} + \mtx{\Phi} - \diag(\vct{\gamma}),
\\
\text{Primal feasibility:~} \widehat{\mtx{X}} \succeq \mtx{0},~\widehat{\mtx{X}} \geq \mtx{0},~\widehat{X}_{ii}\leq 1\text{~for~} 1\leq i \leq n,
\\
\text{Dual feasibility:~} \mtx{\Lambda} \succeq \mtx{0},~\mtx{\Phi}\geq \mtx{0},~\vct{\gamma} \geq \vct{0},
\\
\text{Complimentary slackness:~}\langle \widehat{\mtx{X}}, \mtx{\Lambda} \rangle =0, \langle \widehat{\mtx{X}}, \mtx{\Phi} \rangle =0, \gamma_{i} = 0 \text{~iff~} \widehat{X}_{ii} <1.
\end{cases}
\end{equation}
Notice that $\mtx{\Lambda}$ is determined by $\mtx{\Phi}$ and $\vct{\gamma}$. It remains to construct the primal solution $\widehat{\mtx{X}}$ and dual variables $\mtx{\Phi}$ and $\vct{\gamma}$ for which the above conditions hold.

\paragraph{Primal solution candidate:}
Inspired by numerical simulations, \citet{Cai2014robust} consider a candidate primal solution of the following form:
 \begin{equation}
 \label{eq:primal_candidate}
\widehat{\mtx{X}} = 
\begin{bmatrix} 
\mtx{J}_{n_1} & \mtx{0} & \vct{1}_{n_1}\vct{x}_1^{\top}
\\ 
\mtx{0} & \mtx{J}_{n_2}    &  \vct{1}_{n_2}\vct{x}_2^{\top} 
\\ 
\vct{x}_1\vct{1}_{n_1}^{\top} & \vct{x}_2\vct{1}_{n_2}^{\top} & \vct{x}_1\vct{x}_1^{\top} + \vct{x}_2\vct{x}_2^{\top}
\end{bmatrix}.
\end{equation}
Assume that
\[
-\mtx{A} + \lambda \mtx{J} + (\alpha - \lambda) \mtx{I}
:=
\begin{bmatrix} 
\star & \star & \widetilde{\mtx{Z}}_1
 \\
\star & \star & \widetilde{\mtx{Z}}_2
 \\
\widetilde{\mtx{Z}}_1^{\top}& \widetilde{\mtx{Z}}_2^{\top}& \widetilde{\mtx{W}}
  \end{bmatrix}
\]
Plugging the parameterization \prettyref{eq:primal_candidate} of the primal solution into the convex relaxation \prettyref{eq:SDP2}, we obtain the following quadratic program:
\begin{equation}
\label{eq:cvx2}
\begin{aligned}
&\min_{\vct{x}_1, \vct{x}_2 \in \mathbb{R}^m} && \sum_{i=1}^2\left\langle \vct{x}_i, \widetilde{\mtx{Z}}_i^{\top}\vct{1}_{n_i}\right\rangle+\frac{1}{2}\sum_{i=1}^2 \vct{x}_i^{\top} \widetilde{\mtx{W}} \vct{x}_{i} 
\\
&\text{subject to}       && \vct{x}_1\geq \vct{0}, \quad \vct{x}_2 \geq \vct{0},
\\
&~                                    && \sum_{i=1}^2 \vct{x}_i^{\top}(\vct{e}_j\vct{e}_j^{\top})\vct{x}_i\leq 1, \quad \text{for~} 1\leq j\leq m,
\end{aligned}
\end{equation}
Slightly abusing notation, we use $(\vct{x}_1, \vct{x}_2)$ to denote the solution to the above program. Plugging this pair $(\vct{x}_1, \vct{x}_2)$ into \prettyref{eq:primal_candidate} gives our primal candidate $\widehat{\mtx{X}}$, which by definition automatically satisfies primal feasibility in KKT condition.

As an interesting and useful byproduct, we know that the solutions to the dual problem of \prettyref{eq:cvx2}, denoted by $\vct{\beta}_1, \vct{\beta}_2, \vct{\xi}  \geq \vct{0} \in \mathbb{R}^m$, satisfy the conditions
\[
\begin{cases}
\widetilde{\mtx{W}}\vct{x}_i+\widetilde{\mtx{Z}}_i^{\top}\vct{1}_{l_i}=\vct{\beta}_i-\diag(\vct{\xi})\vct{x}_i, \quad i=1, 2.
\\
\xi_j \left(1- \sum_{i=1}^r\vct{x}_i^{\top}(\vct{e}_j\vct{e}_j^{\top})\vct{x}_i\right)=0, \quad j=1,\ldots, m,
\\
\langle \vct{x}_i, \vct{\beta}_i\rangle=0, \quad i=1,2.
\end{cases}
\]
We then use $\vct{\beta}_1$, $\vct{\beta}_2$ and $\vct{\xi}$ to construct the dual certificates $\mtx{\Phi}$ and $\vct{\gamma}$ as follows
\[
\mtx{\Phi}=
\begin{bmatrix} 
\mtx{0} & \mtx{\Phi}_{12} &  \frac{1}{n_1}\vct{1}_{n_1}\vct{\beta}_1^{\top}
 \\
\mtx{\Phi}_{12}^{\top} & \mtx{0} &\frac{1}{n_2}\vct{1}_{n_2}\vct{\beta}_2^{\top}
 \\
\frac{1}{n_1}\vct{\beta}_1\vct{1}_{n_1}^{\top} & \frac{1}{n_2}\vct{\beta}_2\vct{1}_{n_2}^{\top} & \mtx{0}
\end{bmatrix}, 
\quad
\vct{\gamma} = \begin{bmatrix} \vct{\gamma}_1 \\ \vct{\gamma}_2 \\ \vct{\xi} \end{bmatrix};
\]
see \citet{Cai2014robust} for the expressions of $\vct{\gamma}_1, \vct{\gamma}_2$ and $\mtx{\Phi}_{12}$.

With the primal and dual solutions constructed above, the KKT conditions in~\prettyref{eq:KKT_robust} can be validated with high probability under appropriate  conditions on $p$, $q$, $ n_1 $, $ n_2 $, $m$ and the tuning parameters $\lambda$ and $\alpha$. Doing so proves the main Theorem 3.1 in \cite{Cai2014robust}. A corollary thereof is given below:


\begin{theorem}
	Under the above $(p, q, n_1, n_2, m)$-generalized \SBM, suppose that $n_1=O(n)$ and $n_2=O(n)$, and parametrize $ p,q,m $ as $p= a (\log n) / n$, $q= b (\log n)/n$, and $m = \theta \log n$. If 
	\[
	a - b > C_0(1+\sqrt{\theta})\sqrt{a},
	\]
	then there exist tuning parameters $\lambda$ and $ \alpha $ such that with high probability any solution to \prettyref{eq:SDP2} is of the form \prettyref{eq:output}.
\end{theorem}

This result implies that if $p = O(\log n /n)$ and $q = O(\log n /n)$, the convex relaxation \prettyref{eq:SDP2} is  able to recover the true clusters in presence of $O(\log n)$ outlier nodes. The result can be extended to the more general setting with a growing number of clusters and a general connectivity probability matrix satisfying strong assortativity; see \citet{Cai2014robust} for the details.

Before concluding this section, we note that other notions of robustness have been considered in the literature on convex community detection. For example, one may consider a semirandom model where a monotone adversary is allowed to add edges within communities and delete edges across communities. Weak consistency results have been established under this model \citep{moitra2016robust}. A more general semirandom model has later been introduced by~\cite{makarychev2016learning}; their results on outlier-robustness are comparable to (and sometimes better than) the results in \citet{Cai2014robust}.

%% file: conclusion.tex
\section{Conclusions}

The convex relaxation approaches to community detection have attracted much recent attention. The first and foremost advantage of such approaches is that they provide a \emph{computationally} tractable solution to the worst-case hard combinatorial problems that arise in community detection and clustering. In this survey, we have focused on the \emph{statistical} advantages of the convex relaxation approaches, including the power in providing rigorous and strong recovery guarantees, the robustness against data corruption and model misspecification, and the flexibility in handling a broad range of variations of the community detection problems. We have highlighted several key analytical tools for establishing such properties under certain statistical models of the data. As can be seen, such analysis is largely decoupled from the specific algorithms used to solve the convex programs, which is an important benefit of the convex relaxation framework.

Due to space limit, we did not provide a comprehensive coverage of the fast-growing recent literature in the area of convex community detection. 
Some notable omissions include:
\begin{itemize}
\item Solving the convex programs: This can be done using a variety of existing powerful solvers such as interior point methods and ADMM~\citep{wright_numerical_opt}. By leveraging specific structures (such as sparsity and low-rankness) that arise in community detection, even more efficient solvers have been developed that  scale well to large datasets. For work in this direction we refer to the papers by \cite{CLX2015,JavamardMotanariRicci15,bandeira2016low, Montanri16, Mei17} and the references therein.

%

\item Empirical performance: The convex relaxation approaches have been shown to enjoy competitive performance on various synthetic and real-world benchmarks; see, e.g., the work in  \cite{Cai2014robust, CLX2015,JavamardMotanariRicci15}.
\item Computational barriers in solving problems that are hard on average: The work in \cite{ChenXu14,JavamardMotanariRicci15,HajekWuXu_one_sdp15} studies the interplay between the statistical and computational limits of SDP relaxations. Based on deep but non-rigorous statistical physics arguments, \citet{decelle2011asymptotic} conjectures
that there exists a fundamental computational barrier that prevents all polynomial-time procedures from achieving the  optimal statistical performance. 
Interested readers are referred to the recent surveys in~\citet{moore2017computer, abbe2017community, wu2018statistical} for detailed discussions.


\item Other considerations and extensions: Examples include extracting explicit clustering from the convex relaxation solutions~\citep{Vershynin14,CLX2015}, dealing with an unknown number of clusters and other unknown problem parameters~\citep{YSC2018,CSX2014}, and clustering partially observed or weighted networks~\citep{Jalali2011clustering,lim2017pairwise}.
\end{itemize}


We hope that this survey will stimulate further research exploiting the power of the convex relaxation approaches to network analysis, and in particular provide impetus to the development of efficient implementations of these approaches for practical and large-scale problems.